\documentclass[%
 aip,
 amsmath,amssymb,
 reprint,%
]{revtex4-1}

\usepackage{amsthm}
\usepackage{amssymb}
\usepackage{mathtools}
\usepackage{graphicx}
\usepackage{dcolumn}
\usepackage{bm}

\usepackage[utf8]{inputenc}
\usepackage[T1]{fontenc}
\usepackage{etoolbox}

\makeatletter
\def\@email#1#2{%
 \endgroup
 \patchcmd{\titleblock@produce}
  {\frontmatter@RRAPformat}
  {\frontmatter@RRAPformat{\produce@RRAP{*#1\href{mailto:#2}{#2}}}\frontmatter@RRAPformat}
  {}{}
}%
\makeatother

\newcommand{\m}{\hspace{0.25mm}}

\usepackage{hyperref}
\hypersetup{
  pdftitle={A binomial-adjusted θ-trapezoidal tau-leaping method for chemical kinetic systems},
  colorlinks=true,
  urlcolor=blue,
}

\theoremstyle{plain} 
\newtheorem{theorem}{Theorem}
\newtheorem{definition}{Definition}

\begin{document}

\preprint{}

\title{A binomial-adjusted $\theta$-trapezoidal tau-leaping method\\ for chemical kinetic systems} 

\author{Maximilian Scott}
    \altaffiliation{ms4314@bath.ac.uk}
\author{Sam McCallum}
    \altaffiliation{sm2942@bath.ac.uk}
\author{James Foster}
    \altaffiliation{jmf68@bath.ac.uk}

\affiliation{Department of Mathematical Sciences, University of Bath, UK}

\date{\today}

\begin{abstract}
In \citet{hu2011weaksecond}, the authors propose a weak second order accurate $\theta$-trapezoidal method as an improvement to the classical tau-leaping method. However, for complex chemical kinetic systems their model's underlying analytical assumptions can fail to hold in practice, leading to a significant reduction in its computational effectiveness. Here we introduce a modification to the method which allows it to regain desirable convergence characteristics, both given the original and a relaxed version of the underlying assumptions. We confirm this improvement through a local error analysis and present numerical examples showing it consistently performs either better or comparable to the $\theta$-trapezoidal method. 
\end{abstract}


\maketitle 

\section{Introduction}
\label{sec:section_1}

Generating simulations of chemical kinetic systems plays a critical role in computational cell biology\cite{turner2004stochastic, erdi2014stochastic}, and has increasingly found applications in machine learning\cite{campbell2022continuous, ren2026fast}. The Stochastic Simulation Algorithm (SSA) \cite{gillespie1977exact} is a well-known exact procedure which models each reaction individually at a frequency determined by the Chemical Master Equation (CME). Despite several refinements\cite{cao2005multiscale, vasav2026logtime, liu2005nested}, this method becomes computational inefficient for high-dimensional or stiff systems. 

Alternatively, the tau-leaping method\cite{gillespie2001approximate} simulates many reactions at once by updating the system state only after discrete intervals via Poisson random variables. Nevertheless although tau-leaping may generate simulations faster than SSA for challenging systems, they will no longer be exact and may face issues ranging from numerical instability in stiff systems to negative species counts, directly contradicting modeling assumptions. To alleviate these issues alternative implicit tau-leaping\cite{cao2007adaptive, ahn2013implicit}, binomial\cite{chatterjee2005binomial, marquez2007binomial} or multilevel\cite{moraes2016multilevel} methods have been proposed alongside techniques such as adaptive step-size selection\cite{cao2006efficient} or enforcing positive species counts\cite{cao2007adaptive, anderson2008incorporating} to name a few.

The tau-leaping method is tantamount to the Euler-Maruyama method used to approximate Stochastic Differential Equations (SDE) driven by Brownian motion, in the sense that it offers a weak first order approximation to expected state functionals governed by the CME, i.e. any derived statistics have an error that scales linearly with the size of each jump interval. A natural extension is to derive a higher order method, which is precisely what \citet{hu2011weaksecond} propose with their $\theta$-trapezoidal tau-leaping method. It is a two-step method directly adapted from the one proposed by \citet{anderson2011weak}; the first being a smaller tau-leaping step and the second a carefully constructed trapezoidal step (similar to a second order Runge-Kutta method). The authors go on to show the method is weak second order under a large volume scaling\cite{anderson2011error}, an assumption which essentially states that every species count is above a certain threshold since we would otherwise simulate their associated reactions exactly using SSA. A drawback to this method is that it is roughly twice as computationally expensive as the tau-leaping method, although this should be offset by the lower discretization error.

Inherent in the second trapezoidal step's construction is the possibility that a Poisson random variable has a negative rate, in which case the step is discarded. We call this a degenerate event. If one occurs, the higher order local approximation collapses and propagates globally; the $\theta$-trapezoidal method then performs similarly to the tau-leaping method but at twice the computational cost. Under the large volume scaling, it can be shown that the probability of a degenerate event occurring is small enough that its effect on the error locally is negligible. In experimental settings however, especially in high dimensional, stiff or even oscillatory systems, it is unlikely that the large volume scaling's underlying assumptions to hold throughout.

In this paper we introduce the \textit{Binomial-adjusted $\theta$-Trapezoidal Tau-Leaping} (B$\theta$TTLe) method, a modification of the $\theta$-trapezoidal method allowing a targeted solution to this issue. In the presence of a degenerate event, instead of simply reducing to a tau-leaping method, B$\theta$TTLe utilizes information about its intermediary state to ``backtrack'' over the first step, allowing it to more closely approximate the underlying dynamics. It is constructed in a way that it preserves the global weak second order convergence results under the large volume scaling established for the $\theta$-trapezoidal method while improving on the local error under the traditional scaling. Experimentally, it can be shown to perform as well the $\theta$-trapezoidal method when there are no, or few, degenerate events, and it consistently performs better when significant number occur during a simulation. 

The outline of the paper is as follows. In Section \ref{sec:section_2} we give a brief introduction of the tau-leaping and $\theta$-trapezoidal methods along with the notation that will be used throughout. In Section \ref{sec:section_3} we present our new method alongside a broad derivation and a qualitative explanation of why it works. In Section \ref{sec:section_4} we present local error results comparing the $\theta$-trapezoidal and B$\theta$TTLe methods, and a more quantitative explanation for why it performs better. The proofs for these results can be found in Appendix \ref{sec:appendix_A}. In Section \ref{sec:section_5} we provide several numerical experiments showcasing the performance of the different methods under the presence of degenerate events. 

\section{Tau-leaping and $\theta$-trapezoidal methods}
\label{sec:section_2}

We consider a well-stirred chemical reaction system with $N$ species $\{ S_1, \dots \m, S_N\}$ reacting through $M$ reaction channels $\{ R_1, \dots \m, R_M\}$. Denote the system state at any time $t$ by $\bold{X}_t = \big(X^{\m (1)}_t, \dots \m, X^{\m (N)}_t \big) \in \mathbb{Z}_{\geq 0}^{\m N}$. Following standard kinetic theory arguments for each reaction $R_j$ define the propensity function $a_j: \mathbb{Z}_{\geq 0}^{\m N} \to \mathbb{R}_{\geq 0}$ such that $a_j ( \bold{X}_t ) dt $ is equal to the probability that a reaction $R_j$ occurs in the infinitesimal time interval $[t, t+dt)$ for all times $t$. If this reaction occurs, denote the state change vector by the $\bold{v}_j = \big(v^{\m (1)}_j, \dots \m, v^{\m (N)}_j \big) \in \mathbb{Z}^{\m N}$ such that $v^{\m (i)}_j$ corresponds to the change in the number of $S_i$ molecules if reaction $R_j$ fires, $\forall i = 1, \dots \m, N; \m \m \m j =1, \dots \m, M$. These dynamics can be summarized through the following SDE driven by Poisson processes
\begin{align}
\label{eq:true_SDE}
	\bold{X}_t = \bold{X}_0 + \sum^M_{j=1} \bold{v}_j \mathcal{P}_j \bigg( \int^t_0 a_j ( \bold{X}_s ) ds \bigg),
\end{align}
where each $\mathcal{P}_j$ is a uni-rate Poisson processes. 

Suppose that we wish to simulate a trajectory of $\bold{X}_t$ for $t \in [0,T]$. Define a partition of intervals $[t_n, t_{n+1}]$ such that $\tau := t_{n+1} - t_n$ for all $n=0, \m \dots, \m T \slash \tau - 1$. In this paper, we use the shorthand $\bold{X}_n := \bold{X}_{t_n}$ for all $n$. Note that in general $\tau$ need not be constant, although it is kept so here for simplicity. In order to construct the tau-leaping method, we assume that the propensity functions remain fixed over each interval $[t_n, t_{n+1})$, equivalently
\begin{align*}
	&\mathcal{P}_j \bigg(\int^{t_{n+1}}_0 a_j ( \bold{X}_s ) ds \bigg) -  \mathcal{P}_j \bigg( \int^{t_n}_0 a_j ( \bold{X}_s ) ds \bigg) \\[4pt]
	&\quad \quad \quad \overset{d}{\approx} \mathcal{P}_j \Big( a_j(\bold{X}_{t_n}) \m \tau \Big) 
\end{align*}
where $\overset{d}{\approx}$ means approximately equal in distribution. Combined with \eqref{eq:true_SDE}, it is now straightforward to define each tau-leaping step.
\begin{definition}[Tau-leaping method]
Given an initial state $\bold{X}_0$ and stepsize $\tau$, recursively simulate for all $n=0, \m \dots, \m T \slash \tau - 1$\label{def:tau_leaping}
\begin{align}
	\bold{X}_{n+1} = \bold{X}_n + \sum^M_{j=1} \bold{v}_j \mathcal{P}_j \Big( a_j(\bold{X}_{t_n}) \m \tau \Big),
\end{align}
where the $\mathcal{P}_j \Big( a_j(\bold{X}_{t_n}) \tau \Big)$ are independent Poisson random variables with rate $a_j(\bold{X}_{t_n}) \m \tau$.
\end{definition}
The derivation of the $\theta$-trapezoidal method is analogous to the weak trapezoidal method for SDEs driven by Brownian noise\cite{anderson2011weak}. Let $\theta \in (0,1)$ be a fixed parameter. We then define the $\theta$-trapezoidal method over two steps. The first step is identical to the tau-leaping, but with step size $\theta \m \tau$. The subsequent step naturally has size $(1-\theta) \m \tau$ and utilizes a trapezoidal like approximation of the integral term in \eqref{eq:true_SDE}  in order to compensate for some of the error in the first step, similar to a predictor-corrector method. Note that by construction we require the method to retain the information about the system state at the intermediary step and incorporate it at the endpoint in order to accurately capture the higher order precision.

\begin{definition}[$\theta$-trapezoidal method] Let\label{def:theta_trapezoidal}
\begin{align} \label{eq:alpha_1_2}
	\alpha_1 = \frac{1}{2 \theta}, \quad \quad \text{and} \quad \quad \alpha_2 = \frac{1}{2\theta} - 1 + \theta.
\end{align}
Then given an initial state $\bold{X}_0$ and stepsize $\tau$, recursively simulate for all $n=0, \m \dots, \m T \slash \tau - 1$
\begin{align}
	&\hspace{2.5mm} \bold{X}^*_n \hspace{1mm} = \bold{X}_n + \sum^M_{j=1} \bold{v}_j \mathcal{P}^{\m (1)}_j,\\[3pt]
	&\bold{X}_{n+1} = \bold{X}^*_n +  \sum^M_{j=1} \bold{v}_j \mathcal{P}^{\m (2)}_j,
\end{align}
where the $\mathcal{P}^{\m (1)}_j$ and $\mathcal{P}^{\m (2)}_j$ are independent Poisson random variables such that
\begin{align*}
	&\mathcal{P}^{\m (1)}_j \sim \mathrm{Poisson}\Big( a_j(\bold{X}_n) \theta \tau \Big), \\[3pt]
	&\mathcal{P}^{\m (2)}_j \sim \mathrm{Poisson}\Big( \big( \alpha_1  a_j(\bold{X}^*_n) - \alpha_2  a_j(\bold{X}_n) \big)^+  \tau \Big), 
\end{align*}
and $( \m \cdot \m )^+ = \text{max} \hspace{0.5mm} (0, \m \cdot \m )$.
\end{definition}
A degenerate event occurs whenever the rate of the second Poisson random variable is 0, for any reaction $R_j$, whenever
\begin{align}
\label{eq:degenerate_event_condition}
	\alpha_1 a_j ( \bold{X}^*_n ) - \alpha_2 a_j (\bold{X}_n) \leq 0.
\end{align} 
Both \citet{anderson2011weak}, and \citet{hu2011weaksecond} argue that under certain assumptions, degenerate events occur infrequently enough that their contribution to the method's accuracy can be altogether ignored. Although theoretically robust, these assumptions can fail to hold in practice. Then, it is clear that truncating the rate under a degenerate event in order to keep the method well defined is undesirable; not only do we lose the improved accuracy of a second order method despite the increased computational cost, but we do so in the extreme case that a degenerate event occurs i.e. when a corrector step is most needed. To this date, the optimal way of dealing with degenerate events remains an open problem.

\section{New method} 
\label{sec:section_3}

In order to remedy the issue of degenerate events, an alternative method would need to have an identical first tau-leaping step and somehow retain as closely as possible the form of the second, despite the rate being negative. Rearranging \eqref{eq:degenerate_event_condition}, notice that in the case of a degenerate event 
\begin{align*}
	\frac{1}{2\theta}  a_j ( \bold{X}^*_n ) + \Big( 1 - \frac{1}{2\theta} \Big) a_j ( \bold{X}_n ) \leq \theta a_j ( \bold{X}_n ),
\end{align*}
so define
\begin{align} \label{eq:accept_probability}
	p_j := \frac{\frac{1}{2\theta} a_j ( \bold{X}^*_n ) + \big( 1 - \frac{1}{2\theta} \big) a_j ( \bold{X}_n )}{ \theta a_j ( \bold{X}_n ) },
\end{align}
such that $p_j \in [0, 1]$ for all $\theta \in [1/2, 1)$ given the non-negativity of the propensity function. Next, consider the following identity
\begin{align*}
	\mathrm{Binomial} \big( \mathrm{Poisson} ( a) , b \big) \overset{d}{=} \mathrm{Poisson}( a \times b),
\end{align*}
where $a$ and $b$ are independent. The idea is then to replicate the trapezoidal-like approximation through a binomial distribution where $\mathcal{P}^{\m (1)}_j$ corresponds to the number of trials and $p_j$ the probability of success.  

\begin{definition}[B$\theta$TTLe method]
Let $\alpha_1$, $\alpha_2$ be as defined in \eqref{eq:alpha_1_2}. Then given an initial state $\bold{X}_0$ and stepsize $\tau$, recursively simulate for all $n=0, \m \dots, \m T \slash \tau - 1$ \label{def:bottle}
\begin{align}
	&\hspace{2.5mm} \bold{X}^*_n \hspace{1mm} = \bold{X}_n + \sum^M_{j=1} \bold{v}_j \mathcal{P}^{\m (1)}_j, \\[3pt]
	&\bold{X}_{n+1} = \bold{X}_n +  \sum^M_{j=1} \bold{v}_j \mathcal{Q}_j,
\end{align}
where depending on whether degenerate events have occurred 
\begin{align}
	\mathcal{Q}_j = \left\{ \begin{array}{ccc}  \mathcal{P}^{\m (1)}_j + \mathcal{P}^{\m (2)}_j & \text{if} \quad  \alpha_1  a_j(\bold{X}^*_n) \geq \alpha_2  a_j(\bold{X}_n), \\[3pt] \mathrm{Binomial} \big(  \mathcal{P}^{\m (1)}_j, p_j \big) & \text{otherwise}. \end{array}\right. \nonumber
\end{align}
and the $\mathcal{P}^{\m (1)}_j$ and $\mathcal{P}^{\m (2)}_j$ are identical to those in Definition \ref{def:theta_trapezoidal}.
\end{definition}

The $\theta$-trapezoidal and B$\theta$TTLe methods coincide if there are no degenerate events. In turn this implies, under the same assumptions, that any results and analyses that hold for the former also hold for the latter. Of interest is to what extent they differ when we relax these assumptions.

A degenerate event occurs if a reaction channel's propensity function rapidly decreases after the first tau-leaping step. This is equivalent to many reactions firing in a short space of time, using up a large amount of the underlying reactant species which in turn decreases the propensity function. Intuitively, a negative rate in the second step can be interpreted as a need to backtrack on some of these reactions; in essence canceling them out. This is precisely what the binomial distribution does; the method revisits each reaction and decides whether to accept/reject its occurrence based on a probability measure tailored to how often the reaction channel should have been firing. Nevertheless we cannot directly apply the prior identity linking the binomial and Poisson distributions since $\mathcal{P}^{\m (1)}_j$ and $p_j$ are not independent, yet we show that the error induced by their interaction is small relative to the repercussion of alternatively neglecting the second step.
 
\section{Error analysis}
\label{sec:section_4}

We compare the $\theta$-trapezoidal and B$\theta$TTLe methods under the presence of degenerate events, since they are identical otherwise. To do so, we analyze their Local weak Truncation Error (LTE), or absolute difference in expectation under some functional after any step $[t_n, t_{n+1}]$. This error is benchmarked against the true SDE's \eqref{eq:true_SDE} Markov semigroup expansion with underlying generator\cite{hu2011weakconvergence}
\begin{align}
\label{eq:generator}
	\big( \mathcal{L} \phi \big) (\bold{x}) = \sum^M_{j=1} a_j ( \bold{x} ) \big( \phi ( \bold{x} + \bold{v}_j ) - \phi ( \bold{x} ) \big)
\end{align}
for any functional $\phi: \mathbb{Z}^N \to \mathbb{R}$, i.e. 
\begin{align}
\label{eq:true_SDE_expansion}
\begin{split}
	&\mathbb{E} \big[ \phi( \bold{x}_{t_{n+1}} ) \vert \bold{x}_{t_n} = \bold{x} \big] \\[2pt]
	&\quad  = \phi (\bold{x}) + \tau \big( \mathcal{L} \phi \big) (\bold{x}) + \frac{1}{2} \tau^2 \big( \mathcal{L}^2 \phi \big) (\bold{x}) + O(\tau^3).
\end{split}
\end{align}
As previously mentioned we perform the analysis under the traditional scaling $\tau \to 0$, as opposed to large volume scaling where the system size $V \to \infty$ for the following two reasons. Firstly the convergence analysis would coincide under the large volume scaling since how degenerate events are dealt with is made irrelevant by their theoretically infrequent occurrence. Secondly it is the failure of the underlying assumptions for the large volume scaling which lead to a significant amount of degenerate events occurring. Since we are interested in the cases where they have already occurred, it makes sense to use the traditional scaling to analyze their impact. See \citet{anderson2011error} for an interesting discussion on the two types of scaling. 

The proofs for the following three theorems can be found in the appendix.

\begin{theorem}[LTE of $\theta$-trapezoidal method]
Let $\bold{X}^{\m \theta-\text{trap}}_n$ denote the splitting scheme in Definition \ref{def:theta_trapezoidal} and $\bold{x}_t$ denote the solution of \eqref{eq:true_SDE} over $t \in [t_n, t_{n+1}]$. Assume that initially $\bold{X}^{\m \theta-\text{trap}}_n = \bold{x}_{t_n} = \bold{x}$ and at least one of the reactions is degenerate. Then\label{thm:theorem_1}
\begin{align}
\label{eq:theta_trap_lte}
	\big\vert \mathbb{E}_n \big[ \phi ( \bold{X}^{\m \theta-\text{trap}}_{n+1} ) \big] - \mathbb{E}_n \big[ \phi ( \bold{x}_{t_{n+1}} ) \big] \big\vert \leq C^{\m \theta-\text{trap}} \m \tau 
\end{align}
for any $n$ as $\tau \to 0$, where $C^{\m \theta-\text{trap}} > 0$ is a function of $\theta$, $\phi$ and $\bold{x}$.
\end{theorem}
\begin{theorem}[LTE of B$\theta$TTLe method]
Let $\bold{X}^{\m \text{B$\theta$TTLe}}_n$ denote the splitting scheme in Definition \ref{def:bottle} and $\bold{x}_t$ denote the solution of \eqref{eq:true_SDE} over $t \in [t_n, t_{n+1}]$. Assume that initially $\bold{X}^{\m \text{B$\theta$TTLe}}_n = \bold{x}_{t_n} = \bold{x}$ and at least one of the reactions is degenerate. Then
\begin{align}
\label{eq:bottle_lte}
	\big\vert \mathbb{E}_n \big[ \phi ( \bold{X}^{\m \text{B$\theta$TTLe}}_{n+1} ) \big] - \mathbb{E}_n \big[ \phi ( \bold{x}_{t_{n+1}} ) \big] \big\vert \leq C^{\m \text{B$\theta$TTLe}} \m \tau
\end{align}
for any $n$ as $\tau \to 0$, where $C^{\m \text{B$\theta$TTLe}} > 0$ is a function of $\theta$, $\phi$ and $\bold{x}$.
\end{theorem}
\noindent
\begin{theorem}[LTE uniform weak error bound] For all test functions $\phi$ and $\theta \in [1 / 2, 1)$, assuming that $\bold{X}^{\m \theta-\text{trap}}_n = \bold{X}^{\m \text{new}}_n = \bold{x}_{t_n} = \bold{x}$ and at least one of the reactions is degenerate, the B$\theta$LLe method's LTE upper bound \eqref{eq:bottle_lte} is never greater than the $\theta$-trapezoidal's LTE upper bound \eqref{eq:theta_trap_lte}, i.e. $C^{\m \text{B$\theta$TTLe}} \leq C^{\m \theta-\text{trap}}$ for all $n$, provided the degenerate reaction's reactant species count is large enough.
\end{theorem}

The restriction that the reactant species count is large enough is a lot less prohibitive than in the large volume scaling, where every species count needs to be high. For example, when the propensity function is linear (accounting for the majority of reactions) and $\theta = 1 / 2$ we simply require the reactant species count to be greater than or equal to 2. This is in stark contrast with the large volume scaling which assume that each $X^{\m(i)}_t = O(V)$ where $V \gg 1$. Explicit lower bounds for the reactant species counts can be found in the proof of Theorem 3.

An interesting property of the new method is that we can unwrap it as
\begin{align*}
	&\mathbb{E} \big[ \phi(\bold{X}^{\m \text{new}}_{n+1}) \vert \bold{X}_n = \bold{x} \big] \\[2pt]
	&\quad = \phi( \bold{x} ) + \tau \Big[ \frac{1}{2\theta} \widetilde{\mathcal{L}} \phi + \big(1 - \frac{1}{2\theta} \big) \mathcal{L} \phi \Big] (\bold{x}) + O(\tau^2),
\end{align*}
where
\begin{align*}
\label{eq:shifted_generator}
	\big( \widetilde{\mathcal{L}} \phi \big) (\bold{x}) = \sum^M_{j=1} a_j ( \bold{x} + \bold{v}_j) \big( \phi ( \bold{x} + \bold{v}_j ) - \phi ( \bold{x} ) \big).
\end{align*} 
This is immediate from the proof of Theorem 2. Notice that this shifted version of the generator is similar to $\mathcal{L}$ except it implies that the probability that a reaction $R_j$ occurs in the time interval $[t, t+dt)$ is equal to $a_j (\bold{X}_t + \bold{v}_j) dt$ instead of  $a_j (\bold{X}_t) dt$. Instead of solving the original SDE it offers an approximation to a ``look-forward'' version, reminiscent to an implicit method. In other words, although in the presence of degenerate events we are no longer approximating the correct chemical model, B$\theta$TTLe has the advantage of implementing a correction by momentarily looking ahead of a reaction firing. The fact that it is nevertheless not solving the original SDE is tamed by both the fact that this LTE is assumed \textit{a posteriori} a degenerate event occurring and also damped by the fact that it is only a convex combination of $\mathcal{L}$ and $\widetilde{\mathcal{L}}$. \\

\section{Numerical experiments}
\label{sec:section_5}

The following We test our new method on two examples, and compare them with the tau-leaping and $\theta$-trapezoidal methods. The examples chosen undergo a significant amount of degenerate events; simpler models would be trivial since degenerate events occur infrequently enough that the B$\theta$TTLe and $\theta$-trapezoidal methods coincide.

In order to compare them, we generate multiple ``true'' endpoint samples using SSA\cite{gillespie1977exact} and measure the Fréchet Inception Distance (FID) between the SSA generated $\bold{X}_T^{\m \text{SSA}}$ and those generated by the tau-leaping methods $\bold{X}_T^{\m m}$ . Denote their sample means by $\mu_{\m \text{SSA}}$ and $\mu_m$, and their sample covariance matrices by $\Sigma_{\m \text{SSA}}$ and $\Sigma_m$. Then assuming that we have enough samples we can compare them via 
\begin{align}
	\begin{split}
	&\text{FID} \m (\bold{X}_T^{\m \text{SSA}}, \bold{X}_T^{\m m}) = \| \mu_{\m \text{SSA}} - \mu_{\m m} \|^2_2 \\
	&\quad + \text{Tr} \bigg( \Sigma_{\m \text{SSA}} + \Sigma_{\m m} - 2\Big(\Sigma_{\m \text{SSA}} \m \Sigma_{\m m}\Big)^{\frac{1}{2}} \bigg).
	\end{split}
\end{align}
Essentially, the FID metric allows us to measure simultaneously how close each of methods' sample mean and covariances are to those of the ``true'' endpoints. A more natural metric for measuring global weak error would be of the form
\begin{align*}
	\text{weak error} = \Big\vert \frac{1}{N_{\m \text{SSA}}} \sum^{N_{\m \text{SSA}}}_{j=1} \phi \big( \bold{X}_{T, \m j}^{\m \text{SSA}} \big) -  \frac{1}{N_{\m m}} \sum^{N_{\m m}}_{j=1} \phi \big( \bold{X}_{T, \m j}^{\m m} \big) \Big\vert
\end{align*}
over $N_{\m \text{SSA}}$ samples of the SSA and $N_{\m m}$ of the method, where $\phi$ is any functional of interest. The issue with this metric is that it does not necessarily offer a holistic evaluation of a method's accuracy for any chosen chemical reaction system: we would need to consider a large set of test functions $\phi$ in order to conduct a comprehensive and unbiased analysis. Instead, we can just consider mean and covariances directly since they are both encoded in the FID metric.

Now since our objective is to compare the accuracy of the different methods, we keep all parameters and implementation techniques equal, even if this means foregoing more optimized simulations. For example, we keep step sizes constant and more importantly we do not undertake any significant precautions to counteract the possibility of negative species counts, instead truncating the counts to zero whenever this happens.

Finally, we pay special attention to the computational cost of generating samples using the different methods. Although for any fixed step size $\tau$ the tau-leaping method \ref{def:tau_leaping} may have a larger FID, it achieves it in half the propensity function evaluations as the other two methods. In raw computational cost, it often outperforms them for simpler systems. Thus we compare the FID metric in terms of the number of propensity function evaluations as we vary $\tau$ for a fairer analysis.

The code for the following experiments can be found at (\url{https://github.com/maxscott3/BOTTLe}).

\subsection{Heated monosaccharide-casein system}

This kinetic model describes the way various sugars (glucose, fructose, galactose and tagatose) interact with casein, a protein commonly found in dairy products (specifically with one of its main amino acid's, lysine), through the Maillard reaction\cite{brands2001reactions}. We implement it using the state change vectors and propensity functions derived experimentally in \citet{brands2002kinetic}. This is a model with 11 species and 11 reaction channels.

\begin{figure}[h]
	\includegraphics[scale=0.55]{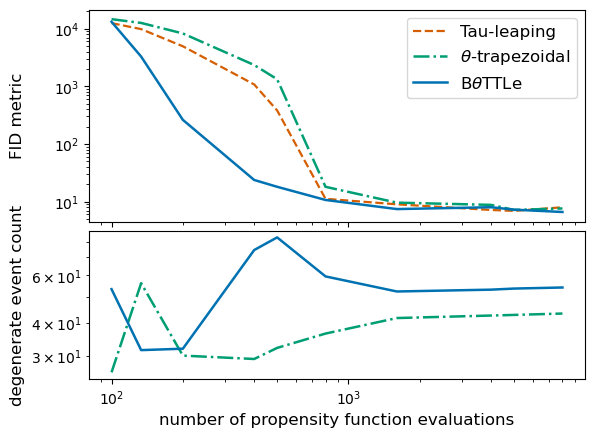}
	\caption{\label{fig:exp1}Experiment comparing FID metric vs propensity function evaluations for methods \ref{def:tau_leaping}, \ref{def:theta_trapezoidal} and \ref{def:bottle} applied to a heated monosaccharide-casein system.}
\end{figure}

We initialized this model with $1.6 \times 10^4 \m \m \text{mmol}$ of glucose and $1.5 \times 10^3 \m \m \text{mmol}$ of lysine residue. We set $T=200$ and measured the FID as we vary $\tau$ from $0.025$ to $2$. We ran $10^3$ simulations of SSA and $10^4$ simulations of the three methods. For the $\theta$-trapezoidal and $B\theta$TTLe methods we set $\theta = 2 \slash 3$ and included an average count of degenerate events in order to confirm that they had indeed occurred.

This kinetic model was chosen since the species count are large enough to satisfy the large volume scaling assumptions used to derive the convergence properties of the $\theta$-trapezoidal method. Nevertheless, although all three methods' endpoints converge to an equilibrium for small $\tau$, as we reduce the number of propensity function evaluations B$\theta$TTle manages to remain numerically stable for longer. The reasoning is that the overall proportion of degenerate event to total reactions increases as we decrease $\tau$, and locally B$\theta$TTLE deviates less from a ``true'' sample whenever they occur by Theorem 3.

\subsection{Polymer growth model}

Here we model the Becker-Döring equations\cite{slemrod2000becker} which, for example, has biophysical applications in the self-assembly of amyloid fibrils\cite{knowles2014amyloid} or prion proteins\cite{masel1999quantifying}. Here we consider a simplified version where mono/polymers attach and detach to form clusters of different lengths. Denote by $S_n$ the polymer of length $n$. Then consider the three following reactions
\begin{align*}
	2 \m S_1 \m &\xrightarrow{\lambda_1} \m S_2, \\[2pt]
	S_1 + S_i \m &\xrightleftharpoons[\lambda_3]{\lambda_2} \m S_{i+1} ,\quad \quad \text{for} \quad i = 1, \m \dots, \m N
\end{align*}
where $N$ is chosen to be the length of the largest polymer. In turn this produces a system with $N$ species and $2N - 1$ reactions channels. 

We consider this model since we can synthetically choose how large of a system we want to analyze. In the following we choose 150 species and 299 reactions, where all the mass is initially concentrated in the monomers, $X_1 = 5000$. We choose reaction rates $\lambda_1 = 10^{-5}$, $\lambda_2 = 0.05$ and $\lambda_3 = 0.01$. We simulate using $T = 0.5$ and vary $\tau =  k \times 10^{-4}$ for $k=1, \dots 10$, using $10^3$ samples of SSA and $10^3$ samples of each of the three methods. We set $\theta=2 \slash 3$ for the $\theta$-trapezoidal and B$\theta$TTLe methods. \\

\begin{figure}[h]
	\includegraphics[scale=0.55]{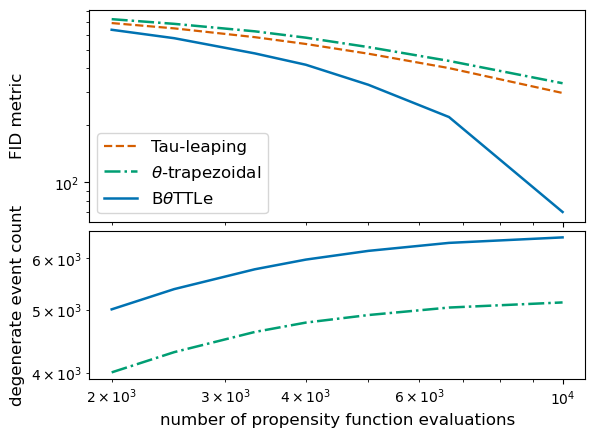}
	\caption{\label{fig:exp2}Experiment comparing FID metric vs propensity function evaluations for methods \ref{def:tau_leaping}, \ref{def:theta_trapezoidal} and \ref{def:bottle} applied to a polymer growth model.}
\end{figure}

The parameters were chosen to make the model intentionally stiff; the samples are relatively slow to transition from the monomers to length 2 polymers, but the reactions speed up once this has occurred. This plays into the explanation of the results; since B$\theta$TTLe approximates a modified system where the reaction is already assumed to have occurred under a degenerate event, it is able to more accurately control the faster reactions. In turn this allows larger choices of $\tau$ to capture both the fast and slow reactions inherent to the system, explaining why B$\theta$TTLe's FID decreases at a faster rate than the tau-leaping and $\theta$-trapezoidal methods.

\section{Conclusion}
\label{sec:section_6}

Overall, the proposed B$\theta$TTLe method is a natural extension and improvement of the $\theta$-trapezoidal method. Analytically it assures identical if not better approximation error rates, regardless of the underlying assumptions on the species' population sizes, suggesting increased efficacy on more complex chemical reaction systems. These theoretical results are backed up by concrete experimental evidence that B$\theta$TTLe is more aptly able to retain desired convergence statistics under the presence of degenerate events. 

An interesting topic of future work is to examine alternative distributions to the binomial in the case of a degenerate event. Indeed, the extraneous $O(\tau)$ term in the LTE of our new method stems solely from the covariance of $\mathcal{P}^{\m (1)}_j$ and $p_j$, an artifact of the moments of a binomial distribution. Additionally, construction of additional predictor-corrector tau-leaping methods over more steps may now be more accessible using a similar methodology as for B$\theta$TTLe; the issue of degenerate events is otherwise magnified by the increased number of steps with potentially negative Poisson rates despite the associated increase in propensity function evaluations.

\appendix
\section{Local truncation error proofs} 
\label{sec:appendix_A}

There is a slight abuse of notation in the statement of Theorems 1 and 2 in that they assume that every single reaction is degenerate. Without loss of generality we can reorder the indexes of the reactions such that only the first $M'$ reactions are degenerate. Then the same result holds once we replace $\mathcal{L}$ with $\mathcal{L'}$ where
\begin{align}
	\big( \mathcal{L'} \phi \big) (\bold{x}) = \sum^{M'}_{j=1} a_j ( \bold{x} ) \big( \phi ( \bold{x} + \bold{v}_j ) - \phi ( \bold{x} ) \big),
\end{align}
and similarly for $\widetilde{\mathcal{L}}$ and $\widetilde{\mathcal{L'}}$. 

The proof strategy is essentially the same for all of the following; we want to expand $\mathbb{E} \big[ \phi( \bold{X}_{n+1} ) \vert \bold{X}_n = \bold{x} \big]$ in terms of $\tau$ and $\bold{x}$ for the splitting scheme, and compare it to the semigroup expansion of the true SDE \eqref{eq:true_SDE}
\begin{align}
\label{eq:true_SDE_expansion_degenerate}
\begin{split}
	\hspace{-2mm} \mathbb{E} \big[ \phi( \bold{x}_{t_{n+1}} ) \vert \bold{x}_{t_n} = \bold{x} \big] = \phi (\bold{x}) + \tau \big( \mathcal{L'} \phi \big) (\bold{x}) + O(\tau^2).
\end{split}
\end{align}
\\
\noindent
\textit{Proof of Theorem 1.} We can apply Taylor's theorem, under the assumption that $\phi$ is sufficiently smooth, so that working backwards
\begin{align*}
	&\mathbb{E} [ \phi ( \bold{X}_{n+1} ) \vert \bold{X}^*_n, \m \bold{X}_n ] \\[2pt]
	&\quad = \phi (\bold{X}^*_n) + \sum^\infty_{p=1} \frac{1}{p!} \phi^{(p)} (\bold{X}^*_n) \Big( \sum^M_{j=M' + 1} \bold{v}_j \mathcal{P}_j^{(2)} \Big)^p \\[3pt]
	& \quad = \phi (\bold{X}^*_n) + \tau \hspace{-3mm} \sum^M_{j=M' + 1} \hspace{-2mm} \Big( \frac{1}{2\theta} a_j(\bold{X}^*_n) + \big(1 - \frac{1}{2\theta} - \theta) \m a_j (\bold{X}_n)  \Big) \\[1pt]
	 &\quad \quad \quad \quad \quad \quad \quad\quad \quad \quad \quad \times \big( \phi ( \bold{x} + \bold{v}_j ) - \phi ( \bold{x} ) \big) \\[1pt]
	 &\quad \quad \quad \quad \quad \m \m+ O(\tau^2),
\end{align*}
given that $\mathbb{E} \big[ \text{Poisson} (\lambda)^n \big] = \lambda + O(\lambda^2)$, $\forall n \geq 2$. Also using the tower property of conditional expectation yields
\begin{align*}
	&\mathbb{E} [ \phi ( \bold{X}_{n+1} ) \vert \bold{X}_n = \bold{x}] \\[2pt]
    &\quad = \mathbb{E} \Big[ \mathbb{E} \big[ \phi ( \bold{X}_{n+1} ) \vert \bold{X}^*_n, \m \bold{X}_n  \big] \m \vert \bold{X}_n = \bold{x} \Big] \\[1pt]
	&\quad = \phi( \bold{x}) + \theta \tau \sum^M_{j=1} a_j ( \bold{x} )  \big( \phi ( \bold{x} + \bold{v}_j ) - \phi ( \bold{x} ) \big) \\[1pt]
	&\quad \quad + (1 - \theta) \m \tau\hspace{-3mm}  \sum^M_{j=M' + 1}\hspace{-2mm}  a_j ( \bold{x} )  \big( \phi ( \bold{x} + \bold{v}_j ) - \phi ( \bold{x} ) \big) + O(\tau^2).
\end{align*}
The result follows by comparing with \eqref{eq:true_SDE_expansion_degenerate} and using the triangle inequality, implying
\begin{align}
\label{eq:constant_lte_theta_trap}
	C^{\m \theta-\text{trap}} := \sum^{M'}_{j=1} \big\vert a_j ( \bold{x} ) \big( \phi ( \bold{x} + \bold{v}_j ) - \phi ( \bold{x} ) \big) \big\vert,
\end{align}
in the limit $\tau \to 0$. \hfill $\Box$ \\
\\
An argument can be made that the local truncation error may be of higher order if we choose $\theta = 1$. However in all but certain special cases, a degenerate event would occur in any reaction channel if the reactant species count decreased over the intermediary step. In turn, the global error would degrade to one no better than the tau-leaping method while requiring twice as many propensity function evaluations.\\
\\
\noindent
\textit{Proof of Theorem 2.} For all $j$, since $\mathcal{P}^{\m (1)}_j$ and $p_j$ are not independent and we require 
\begin{align*}
	\mathbb{E} \Big[ \mathrm{Binomial}\big( \mathcal{P}^{\m (1)}_j, p_j \big)^n \Big], 
\end{align*}
for all $n$, we need to consider how they interact. Firstly by Taylor expanding $a_j$ about $\bold{x}$, we get for all positive integers $m$, $n$
\begin{align*}
	&\mathbb{E}_n \big[ ( \mathcal{P}^{\m (1)}_j)^m  a_j ( \bold{X}^*_n )^n \big] \\[3pt]
	&= \mathbb{E}_n \Big[  ( \mathcal{P}^{\m (1)}_j)^m \sum^n_{k=1} \left( \begin{array}{c} n \\ k \end{array} \right) a_j ( \bold{x} )^k \\[-1pt]
	&\quad \quad \quad \times \bigg( \sum^\infty_{u=1} \frac{1}{u!} a^{\m (u)}_j (\bold{x}) \big( \sum^M_{i=1} \bold{v}_i \mathcal{P}^{\m (1)}_i \big)^u \bigg)^{n-k} \Big] \\[3pt]
	&= \theta \tau a_j ( \bold{x} ) a_j ( \bold{x} + \bold{v}_j )^n + O(\tau^2),
\end{align*}
by matching $\mathcal{P}^{\m (1)}_j$ terms and again using that $\mathbb{E} \big[ \text{Poisson} (\lambda)^n \big] = \lambda + O(\lambda^2)$, $\forall n \geq 2$. Then writing $n^{\m \underbar{k}} = n(n -1)\cdots(n-k+1)$ for the $k$-th falling power of $n$, notice
\begin{align*}
	\mathbb{E}_n \big[ ( \mathcal{P}^{\m (1)}_j)^{\m \underbar{m}}  a_j ( \bold{X}^*_n )^n\big] &= \sum^m_{k=1} s(n,k) \mathbb{E}_n \big[ \Big( ( \mathcal{P}^{\m (1)}_j)^{k}  a_j ( \bold{X}^*_n )^n\Big] \\[1pt]
			&= O(\tau^2),
\end{align*}
for all $m \geq 2$ where $s(n, k)$ are the Stirling numbers of the first kind\cite{wilf2006generatingfunctionology}. As such $\mathbb{E}_n \big[ ( \mathcal{P}^{\m (1)}_j)^{\m \underbar{m}} \m p_j^n\big] = O(\tau^2)$ also for $m \geq 2$. For $m=1$
\begin{align*}
	\mathbb{E}_n \big[ \mathcal{P}^{\m (1)}_j \m p_j \big] = \tau \Big[ \frac{1}{2\theta} a_j (\bold{x} + \bold{v}_j) + \big( 1 - \frac{1}{2\theta} \big) a_j (\bold{x}) \Big] + O(\tau^2).
\end{align*}
Putting everything together
\begin{align*}
	& \mathbb{E} [ \phi ( \bold{X}_{n+1} ) \vert \bold{X}_n = \bold{x}] \\[3pt]
	&= \phi(\bold{x}) + \sum^\infty_{k=1}\frac{1}{k!} \phi^{\m (k)} (\bold{x}) \bigg[ \sum^{M'}_{j=1} \bold{v}_j \mathcal{Q}_j \bigg]^k\\[3pt]
	&= \phi( \bold{x}) + \tau \hspace{-3mm}  \sum^M_{j={M' + 1}} \hspace{-2mm} a_j ( \bold{x} )  \big( \phi ( \bold{x} + \bold{v}_j ) - \phi ( \bold{x} ) \big)  \\[-1pt]
	&\quad + \tau  \sum^{M'}_{j=1}  \Big( \frac{1}{2\theta} a_j( \bold{x} + \bold{v}_j) + \big(1 - \frac{1}{2\theta} \big) a_j ( \bold{x} ) \Big) \\[-1pt]
	&\quad \quad \quad \times \big( \phi ( \bold{x} + \bold{v}_j ) - \phi ( \bold{x} ) \big) \m + O(\tau^2)
\end{align*}
Again the result follows by comparing this expansion with \eqref{eq:true_SDE_expansion_degenerate} in the limit $\tau \to 0$. Here
\begin{align}
\label{eq:constant_lte_new_method}
\begin{split}
	C^{\m \text{B$\theta$TLLe}} := \frac{1}{2\theta} \sum^{M'}_{j=1} &\Big\vert \big( a_j( \bold{x} + \bold{v}_j) - a_j ( \bold{x} ) \big) \\[1pt]
    &\times \big( \phi ( \bold{x} + \bold{v}_j ) - \phi ( \bold{x} ) \big) \Big\vert, 
\end{split}
\end{align}
by the triangle inequality. \hfill $\Box$ \\
\\
\noindent
\textit{Proof of Theorem 3.} Notice that the result follows directly once we are able to show that the induced error for each reaction channel in \eqref{eq:constant_lte_new_method} are bounded by those in \eqref{eq:constant_lte_theta_trap}. This is equivalent to showing that 
\begin{align}
\label{eq:constant_bound_condition}
	\frac{1}{2\theta(1-\theta)} \big\vert a_j (\bold{x} + \bold{v}_j) - a_j (\bold{x}) \big\vert \leq \vert  a_j (\bold{x}) \vert,
\end{align}
for all $j = 1, \dots, M'$. Clearly these only hold given an additional assumption that the reactant species counts are larger than some lower bound, which we argue is small. \\
\\
In the following we solely focus on monomolecular and bimolecular reactions, corresponding to linear and quadratic propensity functions for the reasons discussed in \citet{gillespie1992rigorous}. Essentially, reactions with a larger set of reactant species not only occur very infrequently, but can also be decomposed into more fundamental reactions. Having said this, the following rationale is straightforwardly generalizable; assuming that the propensity function is a polynomial of degree $n$, then the LHS of \eqref{eq:constant_bound_condition} will always be of degree at most $n-1$. Thus there will always exist a lower bound for $\bold{x}$ which, if holds, proves the theorem. Nevertheless, we explicitly state these bounds for the linear and quadratic propensity functions in order to demonstrate that they are small.

Firstly for linear propensity functions, let 
\begin{align*}
	a_j(\bold{x}) := c_j \m x^{\m (m)},
\end{align*}
where $c_j >0$ and $m=1, \dots, N$. Then by definition of the state change vector $\bold{v}_j$, the reactant count decreases by one after the reaction firing, i.e.
\begin{align*}
	a_j(\bold{x} + \bold{v}_j) := c_j \m (x^{\m (m)} - 1).
\end{align*}
Then for linear propensity functions, \eqref{eq:constant_bound_condition} holds provided that 
\begin{align*}
	x^{\m (m)} \geq \frac{1}{2\theta(1 - \theta)}.
\end{align*}
The quadratic propensity functions either have two differing reactant species or require two of the same for the reaction to fire. The former is defined as 
\begin{align*}
	a_j(\bold{x}) := c_j \m x^{\m (m)} \m x^{\m (n)},
\end{align*}
where $m, n=1, \dots, N$. Sufficient conditions for \eqref{eq:constant_bound_condition} to hold are then
\begin{align*}
	x^{\m (m)} \geq \frac{1}{\theta(1 - \theta)}, \quad \text{and} \quad x^{\m (n)} \geq \frac{1}{\theta(1 - \theta)}.
\end{align*}
The other type of quadratic propensity takes the form
\begin{align*}
	a_j(\bold{x}) := \frac{1}{2} c_j \m x^{\m (m)} (x^{\m (m)} -1)
\end{align*}
from which one can show it is sufficient that
\begin{align*}
	x^{\m (m)} \geq \frac{2}{\theta(1 - \theta)},
\end{align*}
taking into account that the state change vector now decreases the reactant count by two. \hfill $\Box$

\section{Discrete diffusion models: image generation}
Recently, tau-leaping methods are increasingly being applied to machine learning models. Here we apply the B$\theta$TTLe method to image generation, using the pretrained MaskGIT model and discrete diffusion implementation of \citet{ren2026fast}. Each $256\times256$ image is represented by a $16\times16$ grid of tokens, taking values in a codebook of size $K=1024$, together with a mask token $m$. The reverse process progressively replaces masks with image tokens. In the one-hot representation, each replacement can be viewed as a reaction with state change $v_{d,k}=e_{d,k}-e_{d,m}$, where $d$ denotes the position and $k$ the destination token.

We retain the model and sampling setup of \citet{ren2026fast}, to which we refer for implementation details (\url{https://github.com/yuchen-zhu-zyc/DiscreteFastSolver}). Specifically, we apply the $\theta$-Trapezoidal and B$\theta$TTLe methods to compute the reverse process of the discrete diffusion reaction system.

Table~\ref{tab:imagenet} compares the $\theta$-trapezoidal and B$\theta$TTLe methods using FID-50k and Inception Score (IS). B$\theta$TTLe achieves a higher IS at every step count considered, with the difference increasing as the number of steps grows. B$\theta$TTLe's FID metric is also comparable or lower than the $\theta$-Trapezoidal FID across all step counts.

\begin{table}[h]
    \centering
    \caption{ImageNet-256 FID-50k and Inception Score (IS) for the MaskGIT model.}
    \label{tab:imagenet}
    \fontsize{12pt}{14pt}\selectfont
    \begin{tabular}{rcccc}
        \hline 
        & \multicolumn{2}{c}{$\theta$-trapezoidal}
        & \multicolumn{2}{c}{B$\theta$TTLe} \\
        \cline{2-5}
        Steps & FID $\downarrow$ & IS $\uparrow$
              & FID $\downarrow$ & IS $\uparrow$ \\
        \hline 
         4 & \textbf{17.1} & 107.1 & 17.3 & \textbf{114.9} \\
         8 & 12.7 & 131.7 & \textbf{10.5} & \textbf{188.4} \\
        16 & 10.1 & 149.0 & \textbf{9.3} & \textbf{225.7} \\
        32 & \textbf{9.1} & 159.0 & \textbf{9.1} & \textbf{235.2} \\
        \hline
    \end{tabular}
\end{table}

\bibliography{references}

\end{document}